\documentclass[journal]{IEEEtran}
\ifCLASSINFOpdf
\else
\fi

\usepackage{graphicx}
\usepackage{amsmath}
\usepackage{amssymb}

\usepackage{mathtools}

\usepackage{color}

\usepackage{url}
 \usepackage{hyperref}
\usepackage{breakurl}

\usepackage{cite}

\usepackage{multirow}
\usepackage{tabularx}
\usepackage{threeparttable}

\usepackage{algorithmic,algorithm}

\makeatletter
\newcommand*{\rom}[1]{\expandafter\@slowromancap\romannumeral #1@}
\makeatother

\newtheorem{theorem}{Theorem}
\newtheorem{remark}{Remark}

\usepackage{arydshln}
\makeatletter
\renewcommand*\env@matrix[1][*\c@MaxMatrixCols c]{%
	\hskip -\arraycolsep
	\let\@ifnextchar\new@ifnextchar
	\array{#1}}
\makeatother

\usepackage{listings}
\makeatletter
\def\lst@makecaption{%
	\def\@captype{table}%
	\@makecaption
}
\makeatother

\begin{document}
%
\title{The Adaptive Dynamic Programming Toolbox}
%
%
%

\author{Xiaowei~Xing
	and~Dong~Eui~Chang
\thanks{The authors are with the School of Electrical Engineering, Korea Advanced Institute of Science and Technology, Daejeon 34141, Republic of Korea (email: xwxing@kaist.ac.kr; dechang@kaist.ac.kr).}
}

\maketitle

\begin{abstract}
The paper develops the Adaptive Dynamic Programming Toolbox (ADPT), which solves optimal control problems for continuous-time nonlinear systems.
Based on the adaptive dynamic programming technique, the ADPT computes optimal feedback controls from the system dynamics in the model-based working mode, or from measurements of trajectories of the system in the model-free working mode without the requirement of knowledge of the system model.
Multiple options are provided  such that the ADPT can accommodate various customized circumstances.
Compared to other popular software toolboxes for optimal control, the ADPT enjoys its computational precision and speed, which is illustrated with its applications to a satellite attitude control problem.
\end{abstract}

\begin{IEEEkeywords}
Optimal control, adaptive dynamic programming, software toolbox.
\end{IEEEkeywords}

%
\IEEEpeerreviewmaketitle

\section{Introduction}
%
%
%
%
\IEEEPARstart{O}{ptimal} control is an important branch in control  engineering.
For continuous-time dynamical systems, finding an optimal feedback control involves solving the so-called Hamilton-Jacobi-Bellman (HJB) equation \cite{Kirk1970}.
For linear systems, the HJB equation becomes the well-known Riccati equation which results in the linear quadratic regulator \cite{LVS2012}.
For nonlinear systems, solving the HJB equation is generally a formidable task due to its inherently nonlinear nature.
As a result, there has been a great deal of research devoted to approximately solving the HJB equation.
Al'brekht proposed a power series method for smooth systems to solve the HJB equation \cite{Albrekht1961}.
With the assumption that the optimal control and the optimal cost function can be represented in Taylor series, by plugging the series expansions of the dynamics, the cost integrand function, the optimal control and the optimal cost function into the HJB equation and collecting terms degree by degree, the Taylor expansions of the optimal control and the optimal cost function can be recursively obtained.
Similar ideas can be found in \cite{GJ1977} and \cite{NSI1971}.
An approach to approximately solve the HJB based on the grid-based discretization of the state and time space is studied in \cite{Peterson1992}.
A recursive algorithm is developed to sequentially improve the control law which converges to the optimal one by starting with an admissible control in \cite{SL1979}.
This recursive algorithm is commonly referred to as policy iteration (PI) and can be also found in \cite{BSW1997,BSW1998,AKL2005}.
The common limitation of these methods is that the complete knowledge of the system is necessary.

In the past few decades, reinforcement learning (RL) \cite{SB1998} has provided a means to design optimal controllers in an adaptive manner from the point of learning.
Adaptive/approximate dynamic programming (ADP), which is an iterative RL-based adaptive optimal control design method, has been proposed in \cite{VL2009,JJ2014,LPC2015}.
An ADP strategy is presented for nonlinear systems with partially unknown dynamics in \cite{VL2009}, and the necessity of the knowledge of system model is fully relaxed in \cite{JJ2014} and \cite{LPC2015}.

Together with the growth of optimal control theory and methods, several software tools for optimal control have been developed.
Notable examples are Nonlinear Systems Toolbox \cite{Krener_NST}, Control Toolbox \cite{GNSB2018}, ACADO \cite{HFD2011}, its successor ACADOS \cite{ACADOS}, and GPSOP-\rom{2} \cite{GPOPS}.
A common feature of these packages is that system equations are used in them.
Besides, optimal controls generated by \cite{GNSB2018,HFD2011,ACADOS,GPOPS} are open-loop such that an optimal control is computed for each initial state.
So, if the initial state changes, optimal controls should be computed again.
In contrast, the Nonlinear Systems Toolbox \cite{Krener_NST} produces an optimal feedback control by solving the HJB equation.

The primary objective of this paper is to develop a MATLAB-based toolbox that solves optimal feedback control problems computationally for nonlinear systems in the continuous-time domain.
More specifically, employing the adaptive dynamic programming technique, we derive a computational methodology to compute approximate optimal feedback controls, based on which we develop the Adaptive Dynamic Programming Toolbox (ADPT).
The ADPT supports two working modes: the model-based mode and the model-free mode.
The knowledge of system equations is required in the model-based working mode.
In the model-free working mode, the ADPT produces the approximate optimal feedback control from measurements of system trajectories, removing the requirement of the knowledge of system equations.
Moreover, multiple options are provided  such that the user can use the toolbox with much flexibility.

The remainder of the paper is organized as follows.
Section \ref{section2} reviews the standard optimal control problem for a class of continuous-time nonlinear systems  and the model-free adaptive dynamic programming technique.
Section \ref{section3} provides implementation details and software features of the ADPT.
In Section \ref{section4}, the ADPT is applied to a satellite attitude control problem in both the model-based working mode and the model-free working mode.
Conclusions and potential future extensions are given in Section \ref{section5}. The codes of the ADPT are available at 
\url{https://github.com/Everglow0214/The_Adaptive_Dynamic_Programming_Toolbox}.

\section{Review of Adaptive Dynamic Programming}
\label{section2}

We review the adaptive dynamic programming (ADP) technique to solve optimal control problems \cite{JJ2017}.
Consider a continuous-time nonlinear system given by
\begin{equation}
	\label{sys}
	\dot{x} = f(x) + g(x)u,
\end{equation}
where $x\in\mathbb{R}^n$ is the state, $u\in\mathbb{R}^m$ is the control, $f:\mathbb{R}^n\to\mathbb{R}^n$ and $g:\mathbb{R}^n\to\mathbb{R}^{n\times m}$ are locally Lipschitz continuous mappings with $f(0)=0$.
It is assumed that the system (\ref{sys}) is stabilizable at $x=0$ in the sense that the system can be locally asymptotically stabilized by a continuous feedback control.
To quantify the performance of a control, an integral cost associated with the system (\ref{sys}) is given by
\begin{equation}
	\label{cost}
	J(x_0,u) = \int_0^\infty(q(x(t))+u(t)^{T}Ru(t))\ \mathrm{d}t,
\end{equation}
where $x_0=x(0)$ is the initial state, $q:\mathbb{R}^n\to\mathbb{R}_{\geq 0}$ is a positive definite function and $R\in\mathbb{R}^{m\times m}$ is a symmetric, positive definite matrix.
A feedback control $u : \mathbb{R}^n \rightarrow \mathbb{R}^m$ is said to be admissible if it stabilizes the system (\ref{sys}) at the origin, and makes the cost $J(x_0,u)$ finite for all  $x_0$ in a neighborhood of $x=0$.

The objective is to find a control policy $u$ that minimizes $J(x_0,u)$ given $x_0$.
Define the optimal cost function $V^{*}:\mathbb{R}^n\to\mathbb{R}$ by
\begin{displaymath}
	V^{*}(x) = \min_{u}J(x,u)
\end{displaymath}
for $x\in\mathbb{R}^n$.
Then, $V^{*}$ satisfies the HJB equation
\begin{displaymath}
	0=\min_{u}\{\nabla V^{*}(x)^{T}(f(x)+g(x)u)+q(x)+u^{T}Ru\},
\end{displaymath}
and the minimizer in the HJB equation is the optimal control which is expressed in terms of $V^{*}$ as
\begin{displaymath}
	u^{*}(x) = -\frac{1}{2}R^{-1}g(x)^T\nabla V^{*}(x).
\end{displaymath}
Moreover, the state feedback $u^{*}$ locally asymptotically stabilizes (\ref{sys}) at the origin and minimizes (\ref{cost}) over all admissible controls \cite{LVS2012}.
Solving the HJB equation analytically is extremely difficult in general except for linear cases.
Hence, approximate or iterative methods are needed to solve the HJB, and a well-known policy iteration algorithm, which is  called policy iteration \cite{SL1979}, is reviewed in Algorithm \ref{algo_PI}.
Let $\{V_i(x)\}_{i\geq 0}$ and $\{u_{i+1}(x)\}_{i\geq 0}$ be the sequences of functions generated by the policy iteration algorithm in Algorithm \ref{algo_PI}.
It is shown in \cite{SL1979} that $V_{i+1}(x)\leq V_i(x)$ for $i\geq 0$, and the limit functions $V(x)=\lim_{i\to\infty}V_i(x)$ and $u(x)=\lim_{i\to\infty}u_i(x)$ are equal to the optimal cost function $V^*$ and the optimal control $u^*$.

\begin{algorithm}[h]
\caption{Policy Iteration}
\label{algo_PI}
\begin{algorithmic}[1]
	\renewcommand{\algorithmicrequire}{\textbf{Input:}}
	\renewcommand{\algorithmicensure}{\textbf{Output:}}
	\newcommand{\algorithmicbreak}{\textbf{break}}
	\newcommand{\BREAK}{\STATE \algorithmicbreak}
	\newcommand{\algorithmiccont}{\textbf{continue}}
	\newcommand{\CONTINUE}{\STATE \algorithmiccont}
	\REQUIRE An initial admissible control $u_0(x)$, and a threshold $\epsilon>0$.
	\ENSURE  The approximate optimal control $u_{i+1}(x)$ and the approximate optimal cost function $V_i(x)$.
	\STATE Set $i \gets 0$.
	\WHILE{$i\geq 0$}
		\STATE Policy evaluation: solve for the continuously differentiable cost function $V_i(x)$ with $V_i(0)=0$ from
		\begin{align}\label{policy_evaluation}
			&\nabla V_i(x)^T(f(x)+g(x)u_i(x)) + q(x) \nonumber\\
			&+u_i(x)^TRu_i(x)=0.
		\end{align}
		\STATE Policy improvement: update the control policy by
		\begin{equation}
			\label{policy_improvement}
			u_{i+1}(x) = -\frac{1}{2}R^{-1}g(x)^T\nabla V_i(x).
		\end{equation}
		\IF{$\Vert u_{i+1}(x)-u_i(x)\Vert\leq\epsilon$ for all  $x$}
		\BREAK
		\ELSE
		\STATE Set $i \gets i+1$.
		\CONTINUE
		\ENDIF
	\ENDWHILE
\end{algorithmic}
\end{algorithm}

Consider approximating the solutions to (\ref{policy_evaluation}) and (\ref{policy_improvement}) by ADP instead of obtaining them exactly.  For this purpose, choose an admissible feedback control $u_0 : \mathbb R^n \rightarrow \mathbb R^m$ for \eqref{sys} and  let $\{V_i(x)\}_{i\geq 0}$ and $\{u_{i+1}(x)\}_{i\geq 0}$ be the sequences of functions generated by the policy iteration algorithm in Algorithm \ref{algo_PI} starting with the  control $u_0(x)$. Choose a bounded time-varying  exploration signal $\eta : \mathbb R \rightarrow \mathbb R^m$, and apply the sum $u_0(x) + \eta(t)$ to \eqref{sys} as follows:
\begin{equation}
	\label{sys_explore_1}
	\dot{x}=f(x)+g(x)(u_0(x)+\eta (t)).
\end{equation}
Assume that solutions to (\ref{sys_explore_1}) are well defined for all positive time.
Let $\mathcal{T} (x,u_0,\eta, [r,s]) = \{ (x(t), u_0(x(t)), \eta(t)) \mid r \leq t \leq s \}$ denote the trajectory $x(t)$ of the system (\ref{sys_explore_1}) with the input $u_0 +\eta$ over the time interval $[r,s]$ with $0\leq r <s$.
The system (\ref{sys_explore_1}) can be rewritten as
\begin{equation}
	\label{sys_explore_2}
	\dot{x}=f(x)+g(x)u_i(x)+g(x)\nu_i (x,t),
\end{equation}
where
\begin{displaymath}
	\nu_i (x,t)=u_0 (x) -u_i (x)+\eta (t).
\end{displaymath}
Combined with (\ref{policy_evaluation}) and (\ref{policy_improvement}), the time derivative of $V_i(x)$ along the trajectory $x(t)$ of (\ref{sys_explore_2}) is obtained as
\begin{equation}
	\label{Vi_dot}
	\dot{V}_i(x) = -q(x) - u_i(x)^{T}Ru_i(x) - 2u_{i+1}(x)^{T}R\nu_i(x,t)
\end{equation}
for $i\geq 0$.
By integrating both sides of (\ref{Vi_dot}) over any time interval $[r,s]$ with $0 \leq r<s$, one gets
\begin{align}
	\label{Vi_integral}
	&V_i(x(s))-V_i(x(r)) \\
	&\quad = -\int_{r}^{s} (q(x)+u_i(x)^{T}Ru_i(x)+2u_{i+1}(x)^{T}R\nu_i (x,t))\ \mathrm{d}\tau. \nonumber
\end{align}

Let $\phi_j:\mathbb{R}^n\to\mathbb{R}$ and $\varphi_j:\mathbb{R}^n\to\mathbb{R}^m$, with $j=1,2,\ldots,$ be two infinite sequences of continuous basis functions on a compact set in $\mathbb{R}^n$ containing the origin as an interior point that vanish at the origin \cite{JJ2017}.
Then, $V_i(x)$ and $u_{i+1}(x)$ for each $i\geq 0$ can be expressed as infinite series of the basis functions.
For each $i\geq 0$ let $\hat{V}_i (x)$ and $\hat{u}_{i+1}(x)$ be approximations of $V_i(x)$ and $u_{i+1}(x)$ given by
\begin{align}
	\hat{V}_i(x) &= \sum_{j=1}^{N_1}c_{i,j}\phi_j(x),
	\label{Vi_app}\\
	\hat{u}_{i+1}(x) &= \sum_{j=1}^{N_2}w_{i,j}\varphi_j(x),
	\label{ui_app}
\end{align}
where $N_1>0$ and $N_2>0$ are integers and $c_{i,j}$, $w_{i,j}$ are coefficients to be found for each $i\geq 0$.
Then, equation (\ref{Vi_integral}) is approximated by $\hat{V}_i (x)$ and $\hat{u}_{i+1}(x)$ as follows:
\begin{align}\label{cw}
	&\sum_{j=1}^{N_1}c_{i,j}(\phi_j(x(s))-\phi_j(x(r))) \nonumber\\
	&+\int_{r}^{s} (2\sum_{j=1}^{N_2}w_{i,j}\varphi_j(x)^{T}R\hat{\nu}_i)\ \mathrm{d}\tau \nonumber\\
	&\quad =- \int_{r}^{s} (q(x)+\hat{u}_i(x)^{T}R\hat{u}_i(x))\ \mathrm{d}\tau,
\end{align}
where
\begin{equation}
	\label{u0:nui}
	\hat{u}_0 = u_0,\quad \hat{\nu}_i = u_0-\hat{u}_i+\eta.
\end{equation}
Suppose that we have available $K$ trajectories $\mathcal {T}(x,u_0,\eta, [r_k,s_k])$, $k = 1, \ldots, K$, where $x(t)$, $u_0(t)$, and $\eta(t)$ satisfy (\ref{sys_explore_2}) over the $K$ time intervals $[r_k, s_k]$, $k = 1, \ldots, K$.
Then, we have $K$ equations of the form (\ref{cw}) for each $i\geq 0$, which can be written as
\begin{equation}
	\label{eqns:ls}
	e_{i,k} = 0, \quad k = 1, \ldots, K,
\end{equation}
where
\begin{align*}
	e_{i,k}\coloneqq &\sum_{j=1}^{N_1}c_{i,j}(\phi_j(x(s_k))-\phi_j(x(r_k)))\\
	&+\int_{r_k}^{s_k} (2\sum_{j=1}^{N_2}w_{i,j}\varphi_j(x)^{T}R\hat{\nu}_i)\ \mathrm{d}\tau\\
	&+ \int_{r_k}^{s_k} (q(x)+\hat{u}_i(x)^{T}R\hat{u}_i(x))\ \mathrm{d}\tau.
\end{align*}
Then, the coefficients $\{c_{i,j}\}_{j = 1}^{N_1}$ and $\{w_{i,j}\}_{j = 1}^{N_2}$ are obtained by minimizing
\begin{displaymath}
	\sum_{k=1}^K \|e_{i,k}\|^2.
\end{displaymath}
In other words, the $K$ equations in (\ref{eqns:ls}) are solved in the least squares sense for the coefficients, $\{c_{i,j}\}_{j = 1}^{N_1}$ and $\{w_{i,j}\}_{j = 1}^{N_2}$.
Thus two sequences $\{\hat{V}_i(x)\}_{i=0}^{\infty}$ and $\{\hat{u}_{i+1}(x)\}_{i=0}^{\infty}$ can be generated from (\ref{cw}).
According to \cite[Cor.~3.2.4]{JJ2017}, for any arbitrary $\epsilon>0$, there exist integers $i^{*}>0$, $N_1^{**}>0$ and $N_2^{**}>0$ such that
\begin{align*}
	\left\Vert\sum_{j=1}^{N_1}c_{i^{*},j}\phi_j(x)-V^{*}(x)\right\Vert &\leq \epsilon,\\
	\left\Vert\sum_{j=1}^{N_2}w_{i^{*},j}\varphi_j(x)-u^{*}(x)\right\Vert &\leq \epsilon
\end{align*}
for all  $x$, if $N_1>N_1^{**}$ and $N_2>N_2^{**}$.

\begin{remark}
	The ADP algorithm relies on only the measurements of states, the initial control policy and the exploration signal, lifting the requirement of knowing the precise system model, while the conventional policy iteration algorithm in Algorithm \ref{algo_PI} requires the knowledge of the exact system model.
	Hence, the ADP algorithm is 100\% data-based and model-free. 
\end{remark}

\begin{remark}
	Equation (\ref{cw}) depends on the initial control $u_0$, the exploration signal $\eta$, the time interval $[r,s]$ as well as the index $i$, where the first three $u_0$, $\eta$, and $[r,s]$ are together equivalent to the trajectory $\mathcal T (x,u_0, \eta, [r,s])$ if the initial state $x(r)$ at $t=r$ is given.
	Hence, we can generate more diverse trajectories by changing $u_0$ and $\eta$ as well as $[r,s]$, and enrich the ADP algorithm accordingly, as follows.
	Suppose that we have available $K$ trajectories  $\mathcal T (x^k,u_0^k, \eta^k, [r_k,s_k])$,  $1 \leq k \leq K$, where $x^k$, $u_0^k$ and $\eta^k$ satisfy (\ref{sys_explore_2}), i.e.,
	\begin{displaymath}
		\dot x^k(t) = f(x^k(t)) + g(x^k(t)) (u_0^k(x^k(t)) + \eta^k(t))
	\end{displaymath}
	for $ r_k \leq t \leq s_k$.
	Then, we have $K$ equations of the form (\ref{cw}) for each $i\geq 0$, which can be written as $e_{i,k} = 0$, $k = 1, \ldots, K$, where
	\begin{align*}
		e_{i,k}\coloneqq &\sum_{j=1}^{N_1}c_{i,j}(\phi_j(x^k(s_k))-\phi_j(x^k(r_k)))\\
		&+\int_{r_k}^{s_k} (2\sum_{j=1}^{N_2}w_{i,j}\varphi_j(x^k)^{T}R\hat{\nu}_i^k)\ \mathrm{d}\tau\\
		&+\int_{r_k}^{s_k} (q(x^k)+\hat{u}_i(x^k)^{T}R\hat{u}_i(x^k))\ \mathrm{d}\tau
	\end{align*}
	with $\hat{u}_0=u_0^k$ and $\hat \nu^k = u_0^k +  \eta^k - \hat u_i$.
	Then, the coefficients $\{c_{i,j}\}_{j=1}^{N_1}$ and $\{w_{i,j}\}_{j=1}^{N_2}$ are obtained by minimizing $\sum_{k=1}^K \|e_{i,k}\|^2$.
	For the sake of simplicity of presentation, however, in this paper we will fix $u_0$ and $\eta$ and vary only the time intervals to generate trajectory data.
\end{remark}

\section{Implementation Details and Software Features}
\label{section3}

We now discuss implementation details and features of the Adaptive Dynamic Programming Toolbox (ADPT).
We provide two modes to generate  approximate optimal feedback controls; one mode requires the knowledge of system model, but the other eliminates this requirement, giving rise to the ADPT's unique capability of handling model-free cases.

\subsection{Implementation of Computational Adaptive Dynamic Programming}

To approximate $V_i(x)$ and $u_{i+1}(x)$ in (\ref{policy_evaluation}) and (\ref{policy_improvement}), monomials composed of state variables are selected as basis functions.
For a pre-fixed number $d\geq 1$, define a column vector $\Phi_d(x)$ by ordering monomials in graded reverse lexicographic order as
\begin{displaymath}
	\Phi_d(x) = (x_1,\ldots,x_n; x_1^2,x_1x_2,\ldots,x_n^2; x_1^3, \ldots, x_n^d)\in\mathbb{R}^{N\times 1},
\end{displaymath}
where $x=(x_1,x_2,\ldots,x_n)\in\mathbb{R}^n$ is the state, $d\geq 1$ is the highest degree of the monomials, and $N$ is given by
\begin{displaymath}
	N=\sum_{i=1}^d{i+n-1 \choose n-1}.
\end{displaymath}
For example, if $n=3$ and $d=3$, the corresponding ordered monomials are
\begin{align*}
	& x_1,x_2,x_3;\\
	& x_1^2,x_1x_2,x_1x_3,x_2^2,x_2x_3,x_3^2;\\
	& x_1^3,x_1^2x_2,x_1^2x_3,x_1x_2^2,x_1x_2x_3,x_1x_3^2,x_2^3,x_2^2x_3,x_2x_3^2,x_3^3.
\end{align*}

According to (\ref{Vi_app}) and (\ref{ui_app}), the cost function $V_i(x)$ and the control $u_{i+1}(x)$ are approximated by $\hat{V}_i(x)$ and $\hat{u}_{i+1}(x)$ which are defined as
\begin{align}
	\hat{V}_i(x) &= c_i \Phi_{d+1}(x),
	\label{vi_hat}\\
	\hat{u}_{i+1}(x) &= W_i \Phi_{d}(x),
	\label{ui_hat}
\end{align}
where $d\geq 1$ is the approximation degree, and $c_i\in\mathbb{R}^{1\times N_1}$ and $W_i\in\mathbb{R}^{m\times N_2}$ are composed of coefficients corresponding to the monomials in $\Phi_{d+1}(x)$ and $\Phi_{d}(x)$ with
\begin{displaymath}
	N_1=\sum_{i=1}^{d+1}{i+n-1 \choose n-1},\quad N_2=\sum_{i=1}^{d}{i+n-1 \choose n-1}.
\end{displaymath}
We take the highest degree of monomials to approximate $V_{i}$   greater by one than the approximation degree since $u_{i+1}$ is obtained by taking the gradient of $V_i$ in (\ref{policy_improvement})  and $g(x)$ is constant in most cases.

\begin{theorem}\label{theorem:main:ADP}
	Let a set of trajectories be defined as $\mathcal{S_T}=\{\mathcal{T} (x,u_0,\eta, [r_k,s_k]), k=1,2,\ldots, K\}$ with  $K\geq1$, and let
	\begin{align*}
		\alpha(x) &= R\eta\Phi_{d}(x)^T,\\
		\beta(x) &= R(u_0(x)+\eta)\Phi_{d}(x)^T,\\
		\gamma(x) &= \Phi_{d-1}(x)\Phi_{d}(x)^T.
	\end{align*}
	Then the coefficients $c_i$ and $W_i$ satisfy
	\begin{equation}
		\label{sys:lin:eqn}
		A_i\begin{bmatrix}
			c_i^T\\
			\mathrm{vec}(W_i)
		\end{bmatrix}= b_i,
	\end{equation}
	where
	\begin{align*}
		A_0 &= \begin{bmatrix}[c:c]
			\Phi_{d+1}^{[r_1,s_1]}(x) & 2\mathrm{vec}(\int_{r_1}^{s_1}\alpha(x)\ \mathrm{d}t)^T\\
			\vdots & \vdots \\
			\Phi_{d+1}^{[r_K,s_K]}(x) & 2\mathrm{vec}(\int_{r_{K}}^{s_{K}}\alpha(x)\ \mathrm{d}t)^T
		\end{bmatrix} \\
		& \in\mathbb{R}^{K\times(N_1+mN_2)},\\
		b_0
		&= \begin{bmatrix}
			-\int_{r_1}^{s_{1}}(q(x)+u_0(x)^TRu_0(x))\ \mathrm{d}t\\
			\vdots\\
			-\int_{r_{K}}^{s_{K}}(q(x)+u_0(x)^TRu_0(x))\ \mathrm{d}t
		\end{bmatrix}\in\mathbb{R}^{K\times 1},
	\end{align*}
	and for $i=1,2,\ldots,$
	\begin{align*}
		A_i &= \begin{bmatrix}[c:c]
			\Phi_{d+1}^{[r_1,s_1]}(x) & 2\mathrm{vec}(\int_{r_1}^{s_{1}}(\beta(x) - RW_{i-1} \gamma(x))\ \mathrm{d}t)^T \\
			\vdots & \vdots \\
			\Phi_{d+1}^{[r_K,s_K]}(x) & 2\mathrm{vec}(\int_{r_K}^{s_{K}}(\beta(x) - RW_{i-1} \gamma(x))\ \mathrm{d}t)^T
		\end{bmatrix} \\
		& \in\mathbb{R}^{K\times(N_1+mN_2)},\\
		b_i &= \begin{bmatrix}
			-\int_{r_1}^{s_{1}}q(x)\ \mathrm{d}t - \langle W_{i-1}^TRW_{i-1},\int_{r_1}^{s_{1}}\gamma(x)\ \mathrm{d}t\rangle\\
			\vdots\\
			-\int_{r_{K}}^{s_{K}}q(x)\ \mathrm{d}t - \langle W_{i-1}^TRW_{i-1},\int_{r_{K}}^{s_{K}}\gamma(x)\ \mathrm{d}t\rangle
		\end{bmatrix} \\
		& \in\mathbb{R}^{K\times 1},
	\end{align*}
	where
	\begin{align*}
		\Phi_{d+1}^{[r_k,s_k]}(x) = \Phi_{d+1}(x(s_k))^T - \Phi_{d+1}(x(r_k))^T
	\end{align*}
	for $k=1,2,\ldots,K$, and the operator $\mathrm{vec}(\cdot)$ is defined as
	\begin{align*}
		\mathrm{vec}(Z) = \begin{bmatrix}
			z_1\\z_2\\\vdots\\z_n
		\end{bmatrix}\in\mathbb{R}^{mn\times 1}
	\end{align*}
	with $z_j\in\mathbb{R}^{m\times 1}$ being the $j$th column of a matrix $Z\in\mathbb{R}^{m\times n}$ for $j=1,\ldots,n$.
\end{theorem}

\begin{IEEEproof}
	See Appendix \ref{appendix_a}.
\end{IEEEproof}

We now give the computational adaptive dynamic programming algorithm in Algorithm \ref{algo_ADP} for practical implementation.
To solve the least squares problem in line 5 in the algorithm, we need to have a sufficiently large number $K$ of trajectories such that the minimization problem can be well solved numerically.
Then the approximate optimal feedback control is generated by the algorithm as $\hat u_{i+1} = W_i \Phi_{d}(x)$.

\begin{algorithm}[h]
	\caption{Computational Adaptive Dynamic Programming}
	\label{algo_ADP}
	\begin{algorithmic}[1]
		\renewcommand{\algorithmicrequire}{\textbf{Input:}}
		\renewcommand{\algorithmicensure}{\textbf{Output:}}
		\newcommand{\algorithmicbreak}{\textbf{break}}
		\newcommand{\BREAK}{\STATE \algorithmicbreak}
		\newcommand{\algorithmiccont}{\textbf{continue}}
		\newcommand{\CONTINUE}{\STATE \algorithmiccont}
		\REQUIRE An approximation degree $d\geq1$, an initial admissible control $u_0(x)$, an exploration signal $\eta(t)$, and a threshold $\epsilon>0$.
		\ENSURE  The approximate optimal control $\hat{u}_{i+1}(x)$ and the approximate optimal cost function $\hat{V}_i(x)$.
		\STATE Apply $u=u_0+\eta$ as the input during a sufficiently long period and collect necessary data.
		\STATE Set $i \gets 0$.
		\WHILE{$i\geq 0$}
		\STATE Generate $A_i$ and $b_i$.
		\STATE Obtain $c_i$ and $W_i$ by solving the minimization problem
		\begin{displaymath}
			\min_{c_i,W_i}\left\Vert A_i
			\begin{bmatrix}
				c_i^T\\
				\mathrm{vec}(W_i)
			\end{bmatrix}-b_i\right\Vert^2.
		\end{displaymath}\\
		\IF{$i\geq1$ and $\left\Vert c_i-c_{i-1} \right\Vert^2 + \left\Vert W_i-W_{i-1} \right\Vert^2\leq\epsilon^2$}
		\BREAK
		\ELSE
		\STATE Set $i \gets i+1$.
		\CONTINUE
		\ENDIF
		\ENDWHILE
		\RETURN $\hat{u}_{i+1}(x)=W_i\Phi_{d}(x)$ and $\hat{V}_i(x)=c_i\Phi_{d+1}(x)$
	\end{algorithmic}
\end{algorithm}

\begin{remark}
	As in the statement of Theorem \ref{theorem:main:ADP}, several integral terms are included in $A_i$ and $b_i$ for $i\geq 0$.
	As in (\ref{u0:nui}), $u_0$ does not get approximated by the basis functions, so the matrices $A_0$ and $b_0$ in Theorem \ref{theorem:main:ADP} are obtained with $x(r_k)$, $x(s_k)$, $\int_{r_k}^{s_{k}}q(x)\ \mathrm{d}t$, $\int_{r_k}^{s_{k}}u_0(x)^TRu_0(x)\ \mathrm{d}t$ and $\int_{r_k}^{s_{k}}\alpha(x)\ \mathrm{d}t$, $1 \leq k \leq K$.
	For $i\geq 1$, the matrices $A_i$ and $b_i$ in Theorem \ref{theorem:main:ADP} need, in addition, $\int_{r_k}^{s_{k}}\beta(x)\ \mathrm{d}t$ and $\int_{r_k}^{s_{k}}\gamma(x)\ \mathrm{d}t$, $1\leq k \leq K$, as well as $W_{i-1}$.
\end{remark}

\begin{remark}
	In the situation where the system dynamic equations are known, the ADPT uses the Runge-Kutta method to simultaneously compute the trajectory points $x(r_k)$ and $x(s_k)$ and the integral terms that appear in $A_i$ and $b_i$.
	In the case when system equations are not known but trajectory data are available, the ADPT applies the trapezoidal method to evaluate these integrals numerically.
	In this case, each trajectory $\mathcal{T}(x,u_0, \eta, [r_k,s_k])$ is represented by a set of its sample points $\{ x(t_{k,\ell}), u_0(t_{k,\ell}), \eta (t_{k,\ell})\}_{\ell = 1}^{L_k}$, where $\{t_{k,\ell}\}_{\ell = 1}^{L_k}$ is a finite sequence that satisfies $r_k = t_{k, 1} < t_{k,2} < \ldots < t_{k, L_k-1} < t_{k, L_k} =  s_k$, and then the trapezoidal method is applied on these sample points to numerically evaluate the integrals over the time interval $[r_k, s_k]$. If intermediate points in the interval $[r_k,s_k]$ are not available so that partitioning the interval $[r_k,s_k]$ is impossible, then we just use the two end points $r_k$ and $s_k$ to evaluate the integral by the trapezoidal method as
\begin{equation}\label{trapezoidal}
\int_{r_k}^{s_k} h(t) dt  \approx \frac{(s_k - r_k) (	h(s_k) + h(r_k))}{2}
\end{equation}
for a function $h(t)$. 
\end{remark}

\subsection{Software Features}
The codes of the ADPT are available at 
\url{https://github.com/Everglow0214/The_Adaptive_Dynamic_Programming_Toolbox}.
\subsubsection{Symbolic Expressions}

It is of great importance for an optimal control package that the user can describe functions such as system equations, cost functions, etc., in a convenient manner.
The idea of the ADPT is to use symbolic expressions.
Consider a optimal control problem, where the system model is in the form \eqref{sys} with
\begin{align}\label{sys:model}
	f(x) =
	\begin{bmatrix}
		x_2\\
		\displaystyle{\frac{-k_1x_1-k_2x_1^3-k_3x_2}{k_4}}
	\end{bmatrix},\quad
	g(x) =
	\begin{bmatrix}
		0\\
		\displaystyle{\frac{1}{k_4}}
	\end{bmatrix},
\end{align}
where $x=(x_1,x_2)\in\mathbb{R}^2$ is the state, $u\in\mathbb{R}$ is the control, and $k_1,k_2,k_3,k_4\in\mathbb{R}$ are system parameters.
The cost function is in the form \eqref{cost} with
\begin{align}\label{cost:l}
	q(x) = 5x_1^2+3x_2^2,\quad R = 2.
\end{align}
Then in the ADPT the system dynamics and the cost function can be defined in lines 1 -- 17 in Listing \ref{list_ex1}.

\begin{lstlisting}[
	caption=An Example of the Model-Based Mode.,
	label=list_ex1,]
  n = 2; % state dimension
  m = 1; % control dimension
  %% Symbolic variables.
  syms x [n,1] real
  syms u [m,1] real
  syms t real
	
  %% Define the system.
  k1 = 3; k2 = 2; k3 = 2; k4 = 5;
  f = [x2;
       (-k1*x1-k2*x1^3-k3*x2)/k4];
  g = [0;
       1/k4];
	
  %% Define the cost function.
  q = 5*x1^2+3*x2^2;
  R = 2;
	
  %% Execute ADP iterations.
  d = 3; % approximation degree
  [w,c] = adpModelBased(f,g,x,n,u,m,q,...
    R,t,d);	
\end{lstlisting}

\subsubsection{Working Modes}

Two working modes are provided in the ADPT; the model-based mode and the model-free mode.
The model-based  mode deals with the situation where the system model is given while the model-free  mode addresses the situation where the system model is not known but only trajectory data are available.
An example of the model-based mode is given in Listing \ref{list_ex1}, where after defining the system model \eqref{sys:model}, the cost function \eqref{cost:l} and the approximation degree $d$ in lines 1 -- 20, the function, $\texttt{adpModelBased}$, returns the coefficients $W_i$ and $c_i$ for the control $\hat u_{i+1}$ and the cost function $\hat V_i$, respectively, in lines 21 -- 22.

An example of the model-free mode is shown in Listing \ref{list_ex2}, where the system model \eqref{sys:model} is assumed to be unknown.
The initial control $u_0$ is in the form of $u_0(x)=-Fx$ with the feedback control gain $F$ defined in line 18.
The exploration signal $\eta$ is composed of four sinusoidal signals as shown in lines 21 -- 22.
A list of two initial states $x(0) = (-3,2)$ and $x(0) = (2.2,3)$ is given in lines 29 -- 30, and a list of the corresponding total time span for simulation is given lines 31 -- 32, where the time interval $[0,6]$ is divided into sub-intervals of size $0.002$ so that trajectory data are recorded every 0.002 sec in lines 38 -- 44.
The time stamps are saved in the column vector $\texttt{t\_save}$ in line 42, and the values of states are saved in the matrix $\texttt{x\_save}$ in line 43, with each row in $\texttt{x\_save}$ corresponding to the same row in $\texttt{t\_save}$.
Similarly, the values of the initial control $u_0$ and the exploration signal $\eta$ are saved in vectors $\texttt{u0\_save}$ and $\texttt{eta\_save}$ in lines 46 -- 47.
These measurements are passed to the function, $\texttt{adpModelFree}$, in lines 51 -- 52 to compute the optimal control and the optimal cost function approximately.

\begin{lstlisting}[
	caption=An Example of the Model-Free Mode.,
	label=list_ex2,
	float,
	floatplacement=t]
  n = 2; % state dimension
  m = 1; % control dimension
  
  %% Define the cost function.
  q = @(x) 5*x(1)^2+3*x(2)^2;
  R = 2;
  
  %% Generate data.
  syms x [n,1] real
  syms t real
  k1 = 3; k2 = 2; k3 = 2; k4 = 5;
  % System dynamics.
  f = [x2;
       (-k1*x1-k2*x1^3-k3*x2)/k4];
  g = [0;
       1/k4];
  
  F = [1, 1] % feedback gain
  
  % Exploration signal.
  eta = 0.8*(sin(7*t)+sin(1.1*t)+...
    sin(sqrt(3)*t)+sin(sqrt(6)*t));
  e = matlabFunction(eta,'Vars',t);
  
  % To be used in the function ode45.
  dx = matlabFunction(f+g*(-F*x+eta),...
    'Vars',{t,x});
  
  xInit = [-3,  2;
           2.2, 3];
  tSpan = [0:0.002:6;
           0:0.002:6];
  odeOpts = odeSet('RelTol',1e-6,...
    'AbsTol',1e-6);
  
  t_save = [];
  x_save = [];
  for i = 1:size(xInit,1)
    [time, states] = ode45(...
      @(t,x)dx(t,x),tSpan(i,:),...
      xInit(i,:),odeOpts);
    t_save = [t_save; time];
    x_save = [x_save; states];
  end
  
  u0_save = -x_save*F;
  eta_save = e(t_save);
  
  %% Execute ADP iterations.
  d = 3; % approximation degree
  [w,c] = adpModelFree(t_save,x_save,n,...
    u0_save,m,eta_save,d,q,R);	
\end{lstlisting}

In the both model-based and model-free modes the approximate control is saved in the file, uAdp.p, that is generated automatically and can be applied by calling $\texttt{u=uAdp(x)}$ without dependence on other files.
Similarly, the user may also check the approximate cost through the file, VAdp.p.

\subsubsection{Options}

Multiple options are provided such that the user may customize optimal control problems in a convenient way.
In the model-based  mode, the user may set option values through the function, $\texttt{adpSetModelBased}$, in a name-value manner before calling $\texttt{adpModelBased}$.
That is, the specified values may be assigned to the named options.
An example is shown in Listing \ref{list_ex3}, where two sets of initial states, time intervals and exploration signals are specified in lines 1 -- 9.
Then, in lines 16 -- 17 the output of $\texttt{adpSetModelBased}$ should be passed to $\texttt{adpModelBased}$ for the options to take effect.
Otherwise, the default values would be used for the options as in lines 21 -- 22 in Listing \ref{list_ex1}.
The options supported by $\texttt{adpSetModelBased}$ are listed in Table \ref{table_adpSetModelBased} in Appendix \ref{appendix_b}.

\begin{lstlisting}[
	caption=A Demonstration of Calling the Function {\textup{\texttt{adpSetModelBased}}}.,
	label=list_ex3,]
  %% The user may specify settings.
  xInit = [-3,  2;
           2.2, 3];
  tSpan = [0, 10;
           0, 8];
  
  syms t real
  eta = [0.8*sin(7*t)+sin(3*t);
           sin(1.1*t)+sin(pi*t)];
  
  adpOpt = adpSetModelBased(...
    'xInit',xInit,'tSpan',tSpan,...
    'explSymb',eta);
  
  %% Execute ADP iterations.
  [w,c] = adpModelBased(f,g,x,n,u,m,q,...
    R,t,d,adpOpt);
\end{lstlisting}

For the command, $\texttt{adpModelFree}$, option values can be modified with the function, $\texttt{adpSetModelFree}$, in the name-value manner.
The options supported by $\texttt{adpSetModelFree}$ are listed in Table \ref{table_adpSetModelFree} in Appendix \ref{appendix_b}.
Among these options, `stride' enables the user to record values of states, initial controls and exploration signals in a high frequency for a long time, while using only a portion of them in the iteration process inside $\texttt{adpModelFree}$.
To illustrate it, let each trajectory in the set $\mathcal{S_T}$ of trajectories in the statement of Theorem \ref{theorem:main:ADP} be represented by two sample points at time $r_k$ and $s_k$, that is, the trapezoidal method evaluates integrals over $[r_k,s_k]$ by taking values at $r_k$ and $s_k$ as in \eqref{trapezoidal}.
Suppose that trajectories in $\mathcal{S_T}$ are consecutive, that is, $s_k=r_{k+1}$ for $k=1,2,\ldots,K-1$.
By setting `stride' to $\delta$ with $\delta=1,2,\ldots,K$, the data used to generate $A_i$ and $b_i$ in Algorithm \ref{algo_ADP} become  $\{\mathcal{T}(x,u_0,\eta,[r_{1+i\delta},s_{(i+1)\delta}]),i\in\mathbb{N},(i+1)\delta\leq K\}$.
For example, consider 3 consecutive trajectories $\mathcal{T}(x,u_0,\eta,[r_k,r_{k+1}])$ with $k=1,2,3$.
If `stride' is set to 1, one will have three equations from \eqref{cw} as follows:
\begin{align*}
	&\sum_{j=1}^{N_1}c_{i,j}(\phi_j(x(r_{k+1})) - \phi_j(x(r_{k}))) \\
	&+ \int_{r_k}^{r_{k+1}} (2\sum_{j=1}^{N_2}w_{i,j}\varphi_j(x)^{T}R\hat{\nu}_i)\ \mathrm{d}\tau \\
	&\quad =- \int_{r_k}^{r_{k+1}} (q(x)+\hat{u}_i(x)^{T}R\hat{u}_i(x))\ \mathrm{d}\tau
\end{align*}
for $k=1,2,3$.
These three equations contribute to three rows of $A_i$ and three rows of $b_i$ as in Theorem \ref{theorem:main:ADP}.
If `stride' is set to 3, then one will have only one equation from \eqref{cw} as follows:
\begin{align}\label{stride:3}
	&\sum_{j=1}^{N_1}c_{i,j}(\phi_j(x(r_{4})) - \phi_j(x(r_{1}))) \nonumber \\
	&+ \int_{r_1}^{r_{4}} (2\sum_{j=1}^{N_2}w_{i,j}\varphi_j(x)^{T}R\hat{\nu}_i)\ \mathrm{d}\tau \nonumber \\
	&\quad =- \int_{r_1}^{r_{4}} (q(x)+\hat{u}_i(x)^{T}R\hat{u}_i(x))\ \mathrm{d}\tau,
\end{align}
where the integrals over $[r_1,r_4]$ are evaluated by the trapezoidal method with the interval $[r_1,r_4]$ partitioned into the three sub-intervals $[r_1,r_2] \cup [r_2, r_3] \cup [r_3, r_4]$, i.e, with the points at $r_1$, $r_2$, $r_3$ and $r_4$.
Equation \eqref{stride:3} will contribute to one row of $A_i$ and one row of $b_i$ as in Theorem \ref{theorem:main:ADP}.
With the assumption that $A_i$ has full rank with `stride' set to 3, by setting `stride' to 3, the number of equations in the minimization problem in Algorithm \ref{algo_ADP} is two less than that with `stride' set to 1, and as a result, the computation load is reduced in the numerical minimization.
It is remarked that with `stride' equal to 3, all the four points at $r_1, \ldots, r_4$ are used by the trapezoidal method to evaluate the integrals over the interval $[r_1, r_4]$ in \eqref{stride:3}, producing a more precise value of integral than the one that would be obtained with the two end points at $r_1$ and $r_4$ only.
An example of calling $\texttt{adpSetModelFree}$ is shown in Listing \ref{list_ex4}.
Similarly, $\texttt{adpModelFree}$ takes the output of $\texttt{adpSetModelFree}$ as an argument to validate the options specified.

\begin{lstlisting}[
	caption=A Demonstration of Calling the Function {\textup{\texttt{adpSetModelFree}}}.,
	label=list_ex4,]
  %% The user may specify settings.
  adpOpt = adpSetModelFree('stride',2);
  
  %% Execute ADP iterations.
  [w,c] = adpModelFree(t_save,x_save,n,...
    u0_save,m,eta_save,d,q,R,adpOpt);	
\end{lstlisting}

\section{Applications to the Satellite Attitude Stabilizing Problem}
\label{section4}

In this section we apply the ADPT to the satellite attitude stabilizing problem because a stabilization problem can be formulated as an optimal control problem.
In the first example, the system model is given and the controller is computed by the function $\texttt{adpModelBased}$.
The same problem is solved again in the second example by the function $\texttt{adpModelFree}$ when the system dynamics is totally unknown.
The source codes for these two examples are available at \url{https://github.com/Everglow0214/The_Adaptive_Dynamic_Programming_Toolbox}.

\subsection{Model-Based Case}
\label{section:model:based}

Let $\mathbb{H}$ denote the set of quaternions and ${\rm S}^3=\{q\in\mathbb{H}\mid \lVert q\rVert=1\}$.
The equations of motion of the continuous-time fully-actuated satellite system are given by
\begin{align}
	\dot{q} &= \frac{1}{2}q\Omega,
	\label{dot_q}\\
	\dot{\Omega} &= \mathbb{I}^{-1}((\mathbb{I}\Omega)\times\Omega)+\mathbb{I}^{-1}u,
	\label{dot_omega}
\end{align}
where $q\in {\rm S}^3$ represents the attitude of the satellite, $\Omega\in\mathbb{R}^3$ is the body angular velocity vector, $\mathbb{I}\in\mathbb{R}^{3\times 3}$ is the moment of inertial matrix and $u\in\mathbb{R}^3$ is the control input.
The quaternion multiplication is carried out for $q\Omega$ on the right-hand side of (\ref{dot_q}) where $\Omega$ is treated as a pure quaternion.
By the stable embedding technique \cite{Chang2018}, the system (\ref{dot_q}) and (\ref{dot_omega}) defined on ${\rm S}^3\times\mathbb{R}^3$ is extended to the Euclidean space $\mathbb{H}\times\mathbb{R}^3$ \cite{Ko2020,Ko2020b} as
\begin{align}
	\dot{q} &= \frac{1}{2}q\Omega-\alpha(| q|^2-1)q,
	\label{dot_q_extend}\\
	\dot{\Omega} &= \mathbb{I}^{-1}((\mathbb{I}\Omega)\times\Omega)+\mathbb{I}^{-1}u,
	\label{dot_omega_extend}
\end{align}
where $q\in\mathbb{H}$, $\Omega\in\mathbb{R}^3$ and $\alpha>0$.

Consider the problem of stabilizing the system (\ref{dot_q_extend}) and (\ref{dot_omega_extend}) at the equilibrium point $(q_e,\Omega_e)=((1,0,0,0),(0,0,0))$.
The error dynamics is given by
\begin{align*}
	\dot{e}_q &= \frac{1}{2}(e_q+q_e)e_{\Omega}-\alpha(| e_q+q_e|^2-1)(e_q+q_e),\\
	\dot{e}_{\Omega} &= \mathbb{I}^{-1}((\mathbb{I}e_{\Omega})\times e_{\Omega}) + \mathbb{I}^{-1}u,
\end{align*}
where $e_q=q-q_e$ and $e_{\Omega}=\Omega-\Omega_e$ are state errors.
Since the problem of designing a stabilizing controller can be solved by designing an optimal controller, we pose an optimal control problem with the cost integral (\ref{cost}) with $q(x)=x^TQx$, where $x=(e_q,e_{\Omega})\in\mathbb{R}^7$ and $Q=2 I_{7\times 7}$, and $R=I_{3\times 3}$.
The inertia matrix $\mathbb{I}$ is set to $\mathbb{I}=\texttt{diag}(0.1029, 0.1263, 0.0292)$.
The parameter $\alpha$ that appears in the above error dynamics is set to $\alpha = 1$.

We set the option `xInit' with three different initial states.
For each initial state, the option `tSpan' is set to $[0,15]$.
We use the option `explSymb' to set exploration signals; see Table \ref{table_adpSetModelBased} in Appendix \ref{appendix_b} for the use of the option `explSysb'.
For the initial control $u_0$, the default initial control is used, which is an LQR controller computed for the linearization of the error dynamics around the origin with the weight matrices $Q=2 I_{7\times 7}$ and $R=I_{3\times 3}$.
We then call the function, $\texttt{adpModelBased}$, to generate controllers of degree $d = 1, 2, 3$.
The computation time taken by the function, $\texttt{adpModelBased}$, to produce the controllers are recorded in Table \ref{table_cost}.
For the purpose of comparison, we also apply Al'brekht's method with the Nonlinear Systems Toolbox (NST) \cite{Krener_NST} to produce controllers of degree $d = 1, 2, 3$ for the same optimal control problem, and record their respective computation time in Table \ref{table_cost}.
For comparison in terms of optimality, we apply the controllers to the system (\ref{dot_q_extend}) and (\ref{dot_omega_extend})  for the initial error state $x_0=((\cos(\theta/2)-1,\sin(\theta/2),0,0),(0,0,0))$ with $\theta=1.99999\pi$ and compute their corresponding values of the cost integral in Table \ref{table_cost}.
Since we do not know the exact optimal value of the cost integral $J(x_0, u)$ for this initial state, we employ the software package called ACADO \cite{HFD2011} to numerically produce the optimal control for this optimal control problem with the given initial state.
We note that both NST and ACADO are model-based.

We can see in Table \ref{table_cost} that ADPT in the model-based mode is superior to NST in terms of optimality, and ADPT (model-based) for $d = 2,3$ is on par with ACADO in terms of optimality.
Notice however that ACADO produces an \textit{open-loop} optimal control for each given initial state, which is an inferior point of ACADO, while ADPT produces a \textit{feedback} optimal control that is independent of initial states.
Moreover, even for the given initial state ACADO takes a tremendous amount of time to compute the open-loop optimal controller.
From these observations, we can say that ADPT in the model-based mode is superior to NST and ACADO in terms of optimality, speed and usefulness all taken into account.

\subsection{Model-Free Case}

Consider solving the same optimal problem as in Section \ref{section:model:based}, but the system dynamics in (\ref{dot_q}) and (\ref{dot_omega}), or equivalently the error dynamics are not available.
Since we do not have real trajectory data available, for the purpose of demonstration we make up some trajectories with four initial states for the error dynamics, where the same initial control $u_0$ and exploration signals $\eta$ are used as the model-based case in Section \ref{section:model:based}.
The simulation for data collection is run over the time interval $[0,20]$ with the recording period being 0.002 sec, producing $10,000 = 20/0.002$ sampled points for each run.
For the function $\texttt{adpModelFree}$, the option of `stride' is set to 4.
Then, the function, $\texttt{adpModelFree}$, is called to generate controllers of degree $d = 1,2,3$, the computation time taken for each of which is recorded in Table \ref{table_cost}.
For the purpose of comparison in terms of optimality, we apply the controllers generated by $\texttt{adpModelFree}$ to the system (\ref{dot_q_extend}) and (\ref{dot_omega_extend}) with the initial error state  $x_0=((\cos(\theta/2)-1,\sin(\theta/2),0,0),(0,0,0))$ with $\theta=1.99999\pi$ and compute the corresponding values of the cost integral; see Table \ref{table_cost} for the values.

From Table \ref{table_cost}, we can see that ADPT in the model-free mode takes more computation time than ADPT in the model-based mode, and the cost integrals by ADPT in the model-free working mode is slightly higher than those in the model-based working mode, since the integrals in the iteration process are evaluated less accurately.
However, ADPT in the model-free mode is superior to NST in terms of optimality and to ACADO in terms of computation time.
More importantly, it is noticeable that the result by model-free ADPT is comparable to model-based ADPT, which shows the power of data-based adaptive dynamic programming and the ADP toolbox.

\begin{figure}[h]
	\begin{center}
		\includegraphics[width=2.5in]{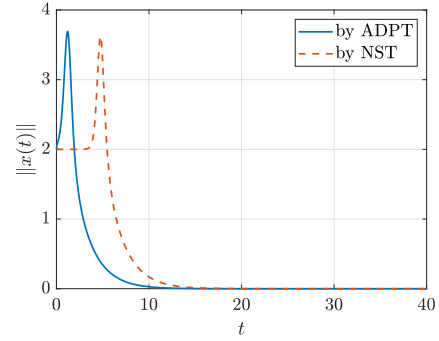}
		\caption{The state errors $\lVert x(t) \rVert$ with the controllers of degree 3 generated by ADPT in the model-free working mode and by NST.}
		\label{fig_err}
	\end{center}
\end{figure}

To see how the computed optimal controller works in terms of stabilization, the norm of the state error under the control with $d = 3$ generated by ADPT in the model-free mode is plotted in Fig. \ref{fig_err} together with the norm of state error by the NST controller with degree 3.
We can see that the convergence to the origin is faster with the model-free ADP controller than with the controller by NST that is model-based.
This comparison result is consistent with the comparison of the two in terms of optimality.


\begin{table}
\begin{center}
	\begin{threeparttable}
		\caption{Costs at $x_0$ and Computation Time by ADPT, NST and ACADO.}
		\label{table_cost}
		\setlength{\tabcolsep}{8pt}
		\setlength{\extrarowheight}{2pt}
		\begin{tabular}{cccc}
			\hline\hline
			& & $J(x_0,u)$\tnote{a} & Time [s]\tnote{b} \\
			\hline
			\multirow{3}{*}{\shortstack{ADPT\\ (model-based)}} & $d=1$ & 37.8259 & 1.6572 \\
			& $d=2$ & 33.6035 & 2.5878 \\
			& $d=3$ & 33.4986 & 11.6869\\
			\hline
			\multirow{3}{*}{\shortstack{ADPT\\ (model-free)}} & $d=1$ & 43.8308 & 0.8923 \\
			& $d=2$ & 36.8319 & 3.5327 \\
			& $d=3$ & 37.4111 & 90.6225 \\
			\hline
			\multirow{3}{*}{NST} & $d=1$ & 208.9259 & 0.2702 \\
			& $d=2$ & 94.6868 & 0.6211 \\
			& $d=3$ & 64.0721 & 3.6201 \\
			\hline
			ACADO & - & 32.6000 & 2359.67 \\
			\hline\hline
		\end{tabular}
		\begin{tablenotes}
			\item[a] $J(x_0,u)$ denotes the integral cost of the corresponding control $u$.
			\item[b] `Time [s]' denotes the computation time taken by the method to obtain the controller.
		\end{tablenotes}
	\end{threeparttable}	
\end{center}
\end{table}

\section{Conclusions and Future Work}
\label{section5}

The Adaptive Dynamic Programming Toolbox, a MATLAB-based package for optimal control for continuous-time nonlinear systems, has been presented.
We propose a computational methodology to approximately produce the optimal control and the optimal cost function by employing the adaptive dynamic programming technique.
The ADPT can work in the model-based mode or in the model-free mode.
The model-based  mode deals with the situation where the system model is given while the model-free  mode handles the situation where the system dynamics are unknown but only system trajectory data are available.
Multiple options are provided for the both modes such that the ADPT can be easily customized.
The optimality, the running speed and the utility of the ADPT are illustrated with a satellite attitude stabilizing problem.

Currently control policies and cost functions are approximated by polynomials in the ADPT.
As mathematical principles of neural networks are being revealed \cite{Gurney1997}, \cite{CC2018}, we plan to use deep neural networks in addition to polynomials in the ADPT to approximately represent optimal controls and optimal cost functions to provide users of the ADPT more options.

\appendices
\section{Proof of Theorem \ref{theorem:main:ADP}}
\label{appendix_a}
\begin{IEEEproof}
	Combining (\ref{cw}), (\ref{vi_hat}) and (\ref{ui_hat}), one has
	\begin{align}\label{new1_i=0}
		&c_0 (\Phi_{d+1}(x(s_{k}))-\Phi_{d+1}(x(r_k))) + 2\int_{r_k}^{s_{k}}\Phi_{d}(x)^{T}W_0^TR\eta\ \mathrm{d}t \nonumber\\
		&\quad = -\int_{r_k}^{s_{k}}(q(x)+u_0(x)^{T}Ru_0(x))\ \mathrm{d}t,
	\end{align}
	and for $i=1,2,\ldots,$
	\begin{align}\label{new1_i>0}
		&c_i (\Phi_{d+1}(x(s_{k}))-\Phi_{d+1}(x(r_k))) \nonumber \\
		&+2\int_{r_k}^{s_{k}}\Phi_{d}(x)^{T}W_i^TR(u_0(x)+\eta)\ \mathrm{d}t \nonumber \\
		&-2\int_{r_k}^{s_{k}}\Phi_{d}(x)^{T}W_i^TRW_{i-1}\Phi_{d}(x)\ \mathrm{d}t \nonumber \\
		&\quad =-\int_{r_k}^{s_{k}}(q(x) + \Phi_{d}(x)^{T}W_{i-1}^TRW_{i-1}\Phi_{d}(x))\ \mathrm{d}t.
	\end{align}
	By applying the property
	\begin{displaymath}
		\langle A, BC\rangle=\langle AC^T, B\rangle=\langle B^TA, C\rangle
	\end{displaymath}
	of the Euclidean inner product defined by $\langle E, F\rangle = \sum_{ij}E_{ij}F_{ij}$ for matrices $E=[E_{ij}]$ and $F= [F_{ij}]$ of equal size, one may rewrite (\ref{new1_i=0}) and (\ref{new1_i>0}) as
	\begin{align}\label{new2_i=0}
		&c_0(\Phi_{d+1}(x(s_{k}))-\Phi_{d+1}(x(r_k))) \nonumber \\
		&+ 2\Big\langle W_0,\, \int_{r_k}^{s_{k}}R\eta\Phi_{d}(x)^T\ \mathrm{d}t\Big\rangle \nonumber\\
		&\quad = -\int_{r_k}^{s_{k}}(q(x)+u_0(x)^{T}Ru_0(x))\ \mathrm{d}t,
	\end{align}
	and for $i=1,2,\ldots,$
	\begin{align}\label{new2_i>0}
		&c_i(\Phi_{d+1}(x(s_{k}))-\Phi_{d+1}(x(r_k))) \nonumber \\
		&+ 2\Big\langle W_i,\, \int_{r_k}^{s_{k}}R(u_0(x)+\eta)\Phi_{d}(x)^T\ \mathrm{d}t\Big\rangle \nonumber \\
		&-2\Big\langle W_i,\, RW_{i-1}\int_{r_k}^{s_{k}}\Phi_{d}(x)\Phi_{d}(x)^T\ \mathrm{d}t\Big\rangle \nonumber \\
		&\quad = - \Big\langle W_{i-1}^TRW_{i-1},\, \int_{r_k}^{s_{k}}\Phi_{d}(x)\Phi_{d}(x)^T\ \mathrm{d}t\Big\rangle \nonumber \\
		&\quad \color{white}= \color{black} -\int_{r_k}^{s_{k}}q(x)\ \mathrm{d}t.
	\end{align}
	Then, the system of linear equations in (\ref{sys:lin:eqn}) readily follows from (\ref{new2_i=0}) and (\ref{new2_i>0}).
\end{IEEEproof}

\section{}
\label{appendix_b}

Here we list options supported by \texttt{adpSetModelBased} and \texttt{adpSetModelFree} in Table \ref{table_adpSetModelBased} and Table \ref{table_adpSetModelFree}, respectively.


\begin{table}[h]
	\begin{center}
		\caption{Options Supported by the Function {\textup{\texttt{adpSetModelBased}.}}}
		\label{table_adpSetModelBased}
		\setlength{\extrarowheight}{2pt}
		\begin{tabular}{p{1.5cm}@{\hskip 1mm}  p{6.5cm}@{\hskip 2mm}}
			\hline\hline
			Name (default) & Description \\
			\hline
			`xInit' \newline ($[\ ]$) & Two default initial states are generated randomly, with the value of each initial state in an open interval $(\text{`xInitMin'},\text{`xInitMax'})$ or $(-\text{`xInitMax'},-\text{`xInitMin'})$. \\
			`xInitNum' \newline ($2$) & The user may choose the number (`xInitNum') and the range (`xInitMin' and `xInitMax') of initial states.\\
			`xInitMin' \newline ($0.3$) & The user may also directly specify initial states (`xInit') by a row vector or a matrix with each row indicating one initial condition.\\
			`xInitMax' \newline ($0.9$) & Trajectories will be calculated starting from these initial states for a time interval (the corresponding row of `tSpan'). The number of rows of `tSpan' automatically matches with `xInitNum' by default.\\
			`tSpan' \newline ($[0\quad 8]$) & Trajectories are calculated by the Runge-Kutta mathod. Settings of the Runge-Kutta mathod (`odeOpt') can be modified by calling the function \texttt{odeSet}.\\
			`odeOpt' & Default value of the option `odeOpt': \texttt{odeset('RelTol',1e-6,'AbsTol',1e-6)}.\\
			\hline
			`u0Symb' \newline ($[\ ]$) & The default initial control is calculated from the linear quadratic method after linearization around the origin.\\
			& The user may also directly specify the initial control (`u0Symb') in a symbolic form.\\
			\hline
			`explSymb' \newline ($[\ ]$) & Default exploration signals are sums of four sinusoidal signals with different frequencies.\\
			`explAmpl' \newline ($0.8$) & The frequencies in the default exploration signals consist of rational numbers and irrational numbers and are chosen randomly from a pre-defined set.\\
			`numFreq' \newline ($4$) & The user may specify the number (`numFreq') and the amplitude (`explAmpl') of sinusoidal signals in one exploration signal.\\
			& The number of default exploration signals automatically matches with the control dimension and the number of initial states.\\
			& All of these default exploration signals are different.\\
			& The user may also directly specify exploration signals (`explSymb') in a symbolic form.\\
			\hline
			`basisOpt' \newline (`mono') & Supported basis functions (`basisOpt'):\newline monomials (`mono').\\
			\hline
			`stride' \newline($1$) & The user may specify data to be used in the iteration process by choosing the stride (`stride') in each trajectory. The number of rows of `stride' automatically matches with `xInitNum' by default.\\
			`crit' \newline ($1$) \newline\newline `epsilon' \newline ($0.001$) & The stop criterion (`crit') can be:\newline $0:\lVert c_i-c_{i-1}\rVert\leq$ `epsilon' (the criterion used in \cite{JJ2017}),\newline $1:\lVert c_i-c_{i-1}\rVert^2+\lVert W_i-W_{i-1}\rVert^2\leq$ `epsilon'$^2$,\newline $2:\lVert c_i-c_{i-1}\rVert\leq$ `epsilon' $\cdot\lVert c_{i-1}\rVert$,\newline $3:\lVert c_i-c_{i-1}\rVert^2+\lVert W_i-W_{i-1}\rVert^2\leq$ `epsilon'$^2\cdot(\lVert c_{i-1}\rVert^2+\lVert W_{i-1}\rVert^2)$.\\
			`maxIter' \newline ($100$) & The maximum number of iterations (`maxIter') can also be specified.\\
			\hline\hline
		\end{tabular}
	\end{center}
\end{table}

\begin{table}
	\begin{center}
		\caption{Options Supported by the Function {\textup{\texttt{adpSetModelFree}.}}}
		\label{table_adpSetModelFree}
		\setlength{\extrarowheight}{2pt}
		\begin{tabular}{p{1.5cm}@{\hskip 1mm}  p{6.5cm}@{\hskip 2mm}}
			\hline\hline
			Name (default) & Description \\
			\hline
			`basisOpt' \newline (`mono') & Supported basis functions (`basisOpt'):\newline monomials (`mono').\\
			\hline
			`stride' \newline($1$) & The user may specify data to be used in the iteration process by choosing the stride (`stride') in each trajectory. The number of rows of `stride' automatically matches with `xInitNum' by default.\\
			`crit' \newline ($1$) \newline\newline `epsilon' \newline ($0.001$) & The stop criterion (`crit') can be:\newline $0:\lVert c_i-c_{i-1}\rVert\leq$ `epsilon' (the criterion used in \cite{JJ2017}),\newline $1:\lVert c_i-c_{i-1}\rVert^2+\lVert W_i-W_{i-1}\rVert^2\leq$ `epsilon'$^2$,\newline $2:\lVert c_i-c_{i-1}\rVert\leq$ `epsilon' $\cdot\lVert c_{i-1}\rVert$,\newline $3:\lVert c_i-c_{i-1}\rVert^2+\lVert W_i-W_{i-1}\rVert^2\leq$ `epsilon'$^2\cdot(\lVert c_{i-1}\rVert^2+\lVert W_{i-1}\rVert^2)$.\\
			`maxIter' \newline ($100$) & The maximum number of iterations (`maxIter') can also be specified.\\
			\hline\hline
		\end{tabular}
	\end{center}
\end{table}

\section*{Acknowledgment}

This work was conducted by Center for Applied Research in Artificial Intelligence(CARAI) grant funded by Defense Acquisition Program Administration(DAPA) and Agency for Defense Development(ADD) (UD190031RD).

\ifCLASSOPTIONcaptionsoff
  \newpage
\fi



%

\bibliographystyle{IEEEtran}
\bibliography{IEEEabrv,ref}


%

%





\end{document}